\documentclass[12pt,twoside,a4paper]{article}
\usepackage{enumerate}
\usepackage{t1enc} 
\usepackage{times}
\usepackage[dvips]{epsfig}
\usepackage{amssymb}
\usepackage{amsmath}
\usepackage{longtable}
\usepackage{makeidx}
\usepackage{euscript}
\usepackage{enumerate}
\usepackage{hyperref}
\DeclareMathAlphabet\scr{U}{scr}{m}{n}
\SetMathAlphabet\scr{bold}{U}{scr}{b}{n}
  \DeclareFontFamily{U}{scr}{\skewchar\font'177}%
  \DeclareFontShape{U}{scr}{m}{n}{<-6>rsfs5<6-8>rsfs7<8->rsfs10}{}%
  \DeclareFontShape{U}{scr}{b}{n}{<-6>rsfs5<6-8>rsfs7<8->rsfs10}{}%

\newcommand{\rr}{\mathbb R}  
\newcommand{\ev}{\mathbb E}

\newcommand{\nn}{\mathbb N}

\newcommand{\pp}{\mathbb P}

\newtheorem{satz}{Theorem}[section]
\newtheorem{prop}[satz]{Proposition}
\newtheorem{theorem}[satz]{Theorem}

\newtheorem{@definition}[satz]{Definition}

\newtheorem{@aufgabe}{Aufgabe}

\newtheorem{@bsp}[satz]{Beispiel}

\newenvironment{mat}{\left (\begin{matrix} } {\end{matrix}\right ) }

\newcommand{\bmat}{\begin{mat}}
\newcommand{\emat}{\end{mat}}
\newcommand{\be}{\begin{enumerate}}
\newcommand{\ee}{\end{enumerate}}
\newcommand{\beq}{\begin{equation}}
\newcommand{\eeq}{\end{equation}}
\newcommand{\bea}{\begin{eqnarray}}
\newcommand{\eea}{\end{eqnarray}}
\newcommand{\beaa}{\begin{eqnarray*}}
\newcommand{\eeaa}{\end{eqnarray*}}

\newcommand{\ep}{\hfill $\Box$}

\renewcommand{\epsilon}{\varepsilon}
\renewcommand{\phi}{\varphi}
\renewcommand{\rho}{\varrho}

\title{{\small\bf LARGE SYSTEMS OF DIFFUSIONS INTERACTING THROUGH THEIR RANKS}\footnote{Research supported in part by NSF grant DMS-0806211.\newline
{\it AMS 2000 subject classifications:} Primary 60J60; secondary 82C22, 65M75.\newline
{\it Keywords and phrases:} Diffusion processes, McKean-Vlasov equation, porous medium equation,
particle method, capital distributions, rank-based market models.}}
\author{{\small BY M. SHKOLNIKOV}\\\quad\\{\small\it Stanford University}}
\date{}

\begin{document}
\maketitle

\begin{abstract}
We study the limiting behaviour of the empirical measure of a system of diffusions interacting through their ranks when the number of diffusions tends to infinity. We prove that the limiting dynamics is given by a McKean-Vlasov evolution equation. Moreover, we show that in a wide range of cases the evolution of the cumulative distribution function under the limiting dynamics is governed by the generalized porous medium equation with convection. The uniqueness theory for the latter is used to establish the uniqueness of solutions of the limiting McKean-Vlasov equation and the law of large numbers for the corresponding systems of interacting diffusions. The implications of the results for rank-based models of capital distributions in financial markets are also explained.   
\end{abstract}

\section{Introduction}

The present article studies the behaviour of the weak solutions to the systems of stochastic differential equations
\begin{eqnarray}\label{pas}
dX_i(t)=\mu(F_{\gamma(t)}(X_i(t)))\;dt+\sigma(F_{\gamma(t)}(X_i(t)))\;dB_i(t),\quad 1\leq i\leq N 
\end{eqnarray}
on an interval $[0,T]$ in the limit $N\rightarrow\infty$. Hereby, $\gamma(t)=\frac{1}{N}\sum_{i=1}^N \delta_{X_i(t)}$ is the empirical measure of the particle system $X_1(t),\dots,X_N(t)$ at time $t$, $F_{\gamma(t)}$ is its cumulative distribution function, $\mu$ and $\sigma$ are measurable functions on $[0,1]$ taking values in $\rr$ and $(0,\infty)$, respectively, and $B_1,\dots,B_N$ are i.i.d. standard Brownian motions. Informally, at any time $t$ the drift and diffusion coefficients of a fixed particle $i$ are determined by its rank in the particle configuration $X_1(t),\dots,X_N(t)$ at time $t$, so that whenever a particle changes its rank, the coefficients change accordingly. The existence and uniqueness of a weak solution to (\ref{pas}) for any $N\in\nn$ was pointed out in \cite{bf} and is essentially due to the results in \cite{bp}, making this description rigorous. For each $N\in\nn$ we fix such a weak solution and denote by $Q^{(N)}$ the probability measure on the space on which it is defined.\\\quad\\
Viewing (\ref{pas}) as an equation for the evolution of the empirical measure $\gamma(t)$ of the particle system on $[0,T]$ we prove that under suitable assumptions on $\mu$, $\sigma$ and the initial positions of the particles the limiting evolution is governed by a McKean-Vlasov evolution equation. Moreover, in the case that the measure in the limiting dynamics is absolutely continuous with respect to the Lebesgue measure on $\rr$ at any time $t\in[0,T]$ we show that the cumulative distribution function of the system evolves according to the generalized porous medium equation with convection:
\begin{eqnarray}\label{pme}
\frac{\partial w}{\partial t}=\frac{d^2}{dx^2} \Sigma(w)-\frac{d}{dx}\Theta(w) 
\end{eqnarray} 
where $\Sigma$ and $\Theta$ are the antiderivatives of $\frac{1}{2}\sigma^2$ and $\mu$, respectively, with $\Sigma(0)=\Theta(0)=0$. For different values of $\Sigma$ and $\Theta$ the equation (\ref{pme}) describes various physical phenomena such as infiltration of water into a porous medium or evaporation of water from soil (see \cite{va} and the references there). Our law of large numbers for the particle systems in (\ref{pas}) shows that the latter can be used to obtain numerical approximations of a continuous $[0,1]$-valued weak solution of the generalized porous medium equation with convection (\ref{pme}) provided that $\Sigma$, $\Theta$ and the initial condition satisfy the assumptions of Theorem 1.2 below. We conjecture that the same is true under more general circumstances and, in particular, for any initial condition which is a cumulative distribution function of a probability measure. This so-called particle method for numerical approximation of solutions of partial differential equations was successfully applied before for the Burgers equation (see \cite{bt1}) and the classical McKean-Vlasov equation as in \cite{mc} (see \cite{bt2}).\\\quad\\
The system of diffusions in (\ref{pas}) and related evolutions of particle systems were studied recently for fixed values of $N\in\nn$ and in some cases for $N=\infty$. They are relevant in the study of capital distributions in financial markets (see e.g. \cite{bf}, \cite{pp}, \cite{cp}, \cite{ip} and \cite{sh2}) and in their discrete time version in the analysis of the Sherrington-Kirkpatrick model of spin glasses (see e.g. \cite{ra}, \cite{aa}, \cite{sh1}). They are also closely related to reflected Brownian motions (see \cite{ip} for the connection) which are widely used as heavy traffic approximations of queueing networks (see e.g. \cite{ha}, \cite{hn}, \cite{tw}, \cite{wi2}).\\\quad\\ 
In the context of capital distributions in financial markets the processes $X_1,\dots,X_N$ stand for logarithmic capitalizations of the firms participating in the market. So, our analysis of the limit $N\rightarrow\infty$ of the described systems gives an understanding of the behaviour of the whole market under the assumption that the number of firms operating in that market is large. It also allows to approximate the evolution of the (logarithmic) capitalization of the $j$-th ranked firm or of the $j$ highest (or lowest) ranked firms in the market for a fixed $1\leq j\leq N$ under the assumption that the number of firms in the market is large.\\\quad\\   
Following McKean's seminal work \cite{mc} systems of diffusing particles in which the drift and diffusion coefficients of each particle are functions of the empirical measure of the whole system and the position of the particle were studied extensively in the context of particle systems with mean field interaction. A good summary of the developments in this direction is given in \cite{ga}. However, in \cite{ga} and the references therein the drift and diffusion coefficients are assumed to be continuous in the position of the particle and in the empirical measure of the system (with respect to the topology of weak convergence on the space of probability measures). This was justified by the continuity of potentials and interaction terms appearing in models of statistical mechanics. In contrast to this, the coefficients in (\ref{pas}) are discontinuous both in the empirical measure $\gamma(t)$ and the position of the particle $X_i(t)$ which necessitates a more delicate analysis. \\\quad\\
To state the main results of the paper we use the following set of notations. For any separable metric space $S$ we denote by $M_1(S)$ the space of probability measures on $S$ endowed with the metric 
\begin{eqnarray}
d(\alpha,\alpha')=\sup_{z:\;\|z\|_\infty+Lip(z)\leq 1}\Big|\int_S z(x)\;\alpha(dx)-\int_S z(x)\;\alpha'(dx)\Big|
\end{eqnarray}
where $\|z\|_\infty$ and $Lip(z)$ denote $\sup_{x\in S}|z(x)|$ and the Lipschitz constant of a function $z:S\rightarrow\rr$, respectively. It is well-known that $d$ metrizes the topology of weak convergence on $M_1(S)$. Moreover, we let $C([0,T],S)$ be the space of continuous functions from $[0,T]$ to $S$, endowed with the topology of uniform convergence. We write $C([0,T])$ for $C([0,T],\rr)$. For $k\in\nn$ we define $C_c^k(\rr)$ as the space of $k$ times continuously differentiable functions on $\rr$ having compact support, set $C_c^\infty(\rr)=\bigcap_{k\geq1} C_c^k(\rr)$ and let $C_c(\rr)$ be the space of continuous functions on $\rr$ with compact support, all equipped with the topology of uniform convergence. For any $t\in[0,T]$ we let $Y_1(t)\leq\dots\leq Y_N(t)$ be the ordered particle system such that $Y_i(t)=X_{\pi_t(i)}(t)$ for some (random) permutation $\pi_t$ of $\{1,\dots,N\}$ depending on $t$. In addition, we introduce the normalized version of the system $X_1,\dots,X_N$ given by $R_i(t)=X_i(t)-Y_{M(N)}(0)$, $1\leq i\leq N$ for $t\in[0,T]$ where $M(N)=1+\frac{N}{2}$ if $N$ is even and $M(N)=\frac{N+1}{2}$ if $N$ is odd. For any fixed $N\in\nn$ we let ${\mathcal R}^{(N)}$ be the distribution of $\frac{1}{N}\sum_{i=1}^N \delta_{R_i(.)}$ as an element of $M_1(C([0,T],M_1(\rr)))$. In the statements on the sequence ${\mathcal R}^{(N)}$, $N\in\nn$ we assume without further mentioning that $\mu$ is strictly decreasing and choose the initial condition of the particles in such a way that $(Y_2(0)-Y_1(0),\dots,Y_N(0)-Y_{N-1}(0))$ is distributed according to the unique invariant distribution of the process $(Y_2(t)-Y_1(t),\dots,Y_N(t)-Y_{N-1}(t))$, $t\in[0,T]$ (see Proposition 1 of \cite{ip} for its existence and uniqueness). Finally, for each $\alpha\in M_1(\rr)$ and $f\in C_c(\rr)$ we write $(\alpha,f)$ for $\int_\rr f\;d\alpha$ and define the operator
\begin{eqnarray}
(L_\alpha f)(x)=f'(x)\mu(F_\alpha(x))+\frac{1}{2}f''(x)\sigma(F_\alpha(x))^2
\end{eqnarray}
acting on $f\in C^2_c(\rr)$ where $F_\alpha$ is the cumulative distribution function of $\alpha$.\\\quad\\
Our main results can be now stated as follows.
{\theorem Let the function $\mu$ be continuously differentiable and such that there exists a constant $\omega_0>0$ with $\mu'(u)\leq -\omega_0$ for all $u\in[0,1]$ and let the function $\sigma^2$ be affine. Then the set $\Xi=\{{\mathcal R}^{(N)},\;N\in\nn\}$ is relatively compact with respect to the topology of weak convergence on $M_1(C([0,T],M_1(\rr)))$ and for any accumulation point ${\mathcal R}^\infty$ of $\Xi$ and any random variable $\rho^\infty$ distributed according to ${\mathcal R}^\infty$ it holds
\begin{eqnarray} \label{mcv2}
(\rho^\infty(t),f)-(\rho^\infty(0),f)=\int_0^t(\rho^\infty(s),L_{\rho^\infty(s)}f)\;ds
\end{eqnarray}
for all $f\in C^3_c(\rr)$ and $t\in[0,T]$ almost surely.} 
{\theorem Suppose that the functions $\mu$ and $\sigma$ are twice continuously differentiable and such that for every solution $\xi$ of the system  
\begin{eqnarray}
&&\forall f\in C_c^\infty(\rr):\quad (\xi(t),f)-(\xi(0),f)=\int_0^t(\xi(s),L_{\xi(s)}f)\;ds, \label{mcv}\\
&&\xi(0)=\lambda \label{mcvic}
\end{eqnarray}
in $C([0,T],M_1(\rr))$ the measures $\xi(t)$, $t\in[0,T]$ are absolutely continuous with respect to the Lebesgue measure on $\rr$ whenever the initial condition $\lambda$ is absolutely continuous with respect to the Lebesgue measure on $\rr$. Then the solution of the system (\ref{mcv}), (\ref{mcvic}) is unique in $C([0,T],M_1(\rr))$ for each such initial condition. If, in addition, the functions $\mu$ and $\sigma^2$ are affine, then the sequence ${\mathcal R}^{(N)}$, $N\in\nn$ converges weakly in $M_1(C([0,T],M_1(\rr)))$ to the Dirac probability measure whose atom is the unique solution of the system (\ref{mcv}), (\ref{mcvic}) in $C([0,T],M_1(\rr))$ with the corresponding initial condition being given in Proposition 3.1 below.\\\quad\\}
{\bf Remark 1.3.\;} A careful reading of the proof shows that the statement of Theorem 1.1 is true for any sequence of particle systems of the form (\ref{pas}) provided that $\mu$ is strictly decreasing and that it holds 
\begin{eqnarray}
&&\sup_{N\in\nn}\frac{1}{N}\sum_{i=1}^N \ev[|Y_i(0)-Y_{M(N)}(0)|]<\infty, \label{cond1}\\
&&\lim_{\epsilon\downarrow0}\limsup_{N\rightarrow\infty}\frac{1}{N^2}\sum_{i,l=1}^N \pp(|Y_i(0)-Y_l(0)|\leq\epsilon)=0 \label{cond2}
\end{eqnarray}
where for any fixed $N\in\nn$ the random vector $(Y_2(0)-Y_1(0),\dots,Y_N(0)-Y_{N-1}(0))$ is distributed according to the unique invariant distribution of the process $(Y_2(t)-Y_1(t),\dots,Y_N(t)-Y_{N-1}(t))$, $t\in[0,T]$. As we show below conditions (\ref{cond1}) and (\ref{cond2}) are satisfied under the assumptions of Theorem 1.1.\\\quad\\
{\bf Remark 1.4.\;} In the case that $\mu$ is twice continuously differentiable and $\sigma$ is a constant function one can deduce the following stochastic representation for an arbitrary solution $\xi\in C([0,T],M_1(\rr))$ of the system (\ref{mcv}), (\ref{mcvic}) whenever the initial condition is absolutely continuous with respect to the Lebesgue measure on $\rr$. Letting $W(t)$, $t\in[0,T]$ be a one-dimensional standard Brownian motion on the time interval $[0,T]$ and $X(t)$, $t\in[0,T]$ be a weak solution to   
\begin{eqnarray}
dX(t)=\mu(F_{\xi(t)}(X(t)))\;dt + \sigma(0)\;dW(t)
\end{eqnarray}
on $[0,T]$ such that the law of $X(0)$ is given by $\xi(0)$ one can proceed as in the proof of Theorem 1.2 below to conclude that for all $t\in[0,T]$ the measure $\xi(t)$ coincides with the law of the random variable $X(t)$. Hence, a standard application of Girsanov's Theorem shows that $\xi(t)$ is absolutely continuous with respect to the Lebesgue measure on $\rr$ for all $t\in[0,T]$. Thus, Theorem 1.2 shows that in this case the solution of the system (\ref{mcv}), (\ref{mcvic}) is unique in $C([0,T],M_1(\rr))$ whenever the initial condition is absolutely continuous with respect to the Lebesgue measure on $\rr$. If, in addition, the function $\mu$ is affine, then the law of large numbers of Theorem 1.2 holds for the particle systems in (\ref{pas}).\\\quad\\ 
Assume now that the functions $\mu$ and $\sigma$ are such that the law of large numbers of Theorem 1.2 applies. In the context of capital distributions in financial markets this means that if the logarithmic capitalizations of the firms in the market follow the dynamics in (\ref{pas}) and the number of the firms is large, then the evolution of the empirical measure $\frac{1}{N}\sum_{i=1}^N \delta_{R_i(t)}$ of normalized logarithmic capitalizations is approximately given by the unique solution of the system (\ref{mcv}), (\ref{mcvic}) with the initial condition of Proposition 3.1 below. Moreover, the evolution of the capitalization of the $j$-th ranked firm or of the $j$ highest (or lowest) ranked firms in the market can be approximated by the evolution of the $\frac{N-j+1}{N}$ quantile or the $\frac{N}{N},\dots,\frac{N-j+1}{N}$ quantiles (or $\frac{1}{N},\dots,\frac{j}{N}$ quantiles) of the solution of the system (\ref{mcv}), (\ref{mcvic}) with the initial condition of Proposition 3.1, respectively. In addition, the proof of Theorem 1.2 shows that the cumulative distribution function of the probability measure $\frac{1}{N}\sum_{i=1}^N \delta_{R_i(t)}$, which describes the fraction of firms whose capitalizations are below a certain threshold, can be approximately described by the unique generalized solution to the Cauchy problem for the generalized porous medium equation with convection (\ref{pme}) in the sense of \cite{gi}, with the initial condition being given by the cumulative distribution function of the measure $\lambda^\infty$ in Proposition 3.1. \\\quad\\
The rest of the paper is organized as follows. In section 2 we give the proof of Theorem 1.1 which relies on a characterization of compact subsets of $C([0,T],M_1(\rr))$ obtained in \cite{ga}, a characterization of tight sequences of probability measures on $C([0,T])$ as in \cite{sv} and results on convergence of semimartingales in the spirit of \cite{js}. In section 3 we determine the appropriate initial condition for the limiting dynamics using Lindeberg's Central Limit Theorem in Proposition 3.1 and present the proof of Theorem 1.2 subsequently. In the latter we use a computation similar to the one in section 1 of \cite{jo} together with the analysis of the Fokker-Planck equation in \cite{ga} to demonstrate that under the limiting dynamics the cumulative distribution function evolves according to the generalized porous medium equation with convection (\ref{pme}). Using the results of \cite{gi} on the latter we obtain the uniqueness of solutions for the system (\ref{mcv}), (\ref{mcvic}) for any absolutely continuous initial condition and as a consequence the law of large numbers of Theorem 1.2. 
 
\section{Proof of Theorem 1.1}

Before proving Theorem 1.1 we recall the results of \cite{ip} on the invariant distribution of the gap process $(Y_2(t)-Y_1(t),\dots,Y_N(t)-Y_{N-1}(t))$, $t\in[0,T]$ for a fixed $N\in\nn$. We remark that the results of \cite{ip} are applicable here, since by subtracting $\frac{1}{N}\sum_{i=1}^N \mu\Big(\frac{i}{N}\Big)t$ from $X_1(t),\dots,X_N(t)$ for all $t\in[0,T]$ and by reversing the order of the labels of the particles we can transform the particle system in (\ref{pas}) into an instance of the particle systems considered in \cite{ip}. It is shown in Proposition 1 of \cite{ip} that if $\mu$ is strictly decreasing, an invariant distribution exists and is unique. Moreover, if the function $\sigma^2$ is affine, then under the invariant distribution the joint law of $Y_2(0)-Y_1(0),\dots,Y_N(0)-Y_{N-1}(0)$ is that of independent exponential random variables with parameters 
\begin{eqnarray*}
a^{(N)}_i=\frac{4i(N-i)}{N}\cdot\frac{\frac{1}{i}\sum_{j=1}^i \mu\Big(\frac{j}{N}\Big)-\frac{1}{N-i}\sum_{j=i+1}^N \mu\Big(\frac{j}{N}\Big)}
{\sigma\Big(\frac{i}{N}\Big)^2+\sigma\Big(\frac{i+1}{N}\Big)^2},\quad 1\leq i\leq N-1
\end{eqnarray*} 
by Proposition 3 of \cite{ip}. If, in addition, $\mu$ is as in Theorem 1.1, then a straightforward computation using the inequalities $\mu(u)-\mu\Big(\frac{i}{N}\Big)\geq \omega_0\Big(\frac{i}{N}-u\Big)$ for all $u\in\Big[0,\frac{i}{N}\Big]$ and $\mu\Big(\frac{i}{N}\Big)-\mu(u)\geq \omega_0\Big(u-\frac{i}{N}\Big)$ for all $u\in\Big[\frac{i}{N},1\Big]$ shows that
\begin{eqnarray} \label{plb}
4\cdot\frac{\frac{1}{i}\sum_{j=1}^i \mu\Big(\frac{j}{N}\Big)-\frac{1}{N-i}\sum_{j=i+1}^N \mu\Big(\frac{j}{N}\Big)}
{\sigma\Big(\frac{i}{N}\Big)^2+\sigma\Big(\frac{i+1}{N}\Big)^2}\geq\frac{\omega_0}{\sup_{u\in[0,1]}\sigma(u)^2}
\end{eqnarray}
for all $1\leq i\leq N-1$, $N\in\nn$.\\\quad\\
{\it Proof of Theorem 1.1.} 1) By Prokhorov's Theorem it suffices to prove that the sequence ${\mathcal R}^{(N)}$, $N\in\nn$ is tight to show that $\Xi$ is relatively compact. To this end, we fix an arbitrary $\epsilon>0$ and a countable dense subset $\{f_1,f_2,\dots\}$ of $C_c(\rr)$ contained in $C^2_c(\rr)$. From the proof of Lemma 1.3 in \cite{ga} we see that it is enough to find a compact set $K_0$ in $M_1(\rr)$ and compact sets $K_1,K_2,\dots$ in $C([0,T])$ such that for all $N\in\nn$:
\begin{eqnarray}
&&{\mathcal R}^{(N)}(\{\xi\in C([0,T],M_1(\rr))|\forall t\in[0,T]:\;\xi(t)\in K_0\})\geq 1-\epsilon,\\
&&{\mathcal R}^{(N)}(\{\xi\in C([0,T],M_1(\rr))|(\xi(.),f_r)\in K_r\})\geq 1-\epsilon\cdot 2^{-r},\; r\geq1. 
\end{eqnarray}  
To define $K_0$ we introduce the function $\phi(x)=\sqrt{1+x^2}$ and apply Ito's formula to compute
\begin{eqnarray*}
&&d(\rho^{(N)}(t),\phi)=\frac{1}{N}\sum_{i=1}^N \phi'(R_i(t))\sigma(F_{\rho^{(N)}(t)}(R_i(t))) dB_i(t)\\
&&+\frac{1}{N}\sum_{i=1}^N \Big(\phi'(R_i(t))\mu(F_{\rho^{(N)}(t)}(R_i(t)))+\frac{1}{2}\phi''(R_i(t))\sigma(F_{\rho^{(N)}(t)}(R_i(t)))^2\Big) dt
\end{eqnarray*}  
where $\rho^{(N)}(t)=\frac{1}{N}\sum_{i=1}^N \delta_{R_i(t)}$. The boundedness of $\phi'$, $\phi''$, $\mu$ and $\sigma$ shows that there exists a constant $C_1>0$ such that
\begin{eqnarray*}
\frac{1}{N}\sum_{i=1}^N \Big(\phi'(R_i(t))\mu(F_{\rho^{(N)}(t)}(R_i(t)))+\frac{1}{2}\phi''(R_i(t))\sigma(F_{\rho^{(N)}(t)}(R_i(t)))^2\Big)\leq C_1
\end{eqnarray*}
for all $t\in[0,T]$ and $N\in\nn$. Moreover, for any fixed $N\in\nn$ the process
\begin{eqnarray}
Z(t)=\frac{1}{N}\sum_{i=1}^N \int_0^t\phi'(R_i(s))\sigma(F_{\rho^{(N)}(s)}(R_i(s))) dB_i(s),\quad t\in[0,T] 
\end{eqnarray}
is a continuous martingale. Applying Doob's maximal inequality for non-negative continuous submartingales and Jensen's inequality we obtain for all $A>0$:
\begin{eqnarray*}
&&Q^{(N)}\Big(\sup_{t\in[0,T]} (\rho^{(N)}(t),\phi)\geq A+C_1 T\Big)\leq Q^{(N)}\Big(\sup_{t\in[0,T]} (b^{(N)}+|Z(t)|)\geq A\Big)\\
&&\leq\frac{1}{A}\ev^{Q^{(N)}}[b^{(N)}]+\frac{1}{A}\ev^{Q^{(N)}}\Big[Z(T)^2\Big]^{\frac{1}{2}}
\end{eqnarray*}
where $b^{(N)}=(\rho^{(N)}(0),\phi)$. Provided that we can show that $\ev^{Q^{(N)}}[b^{(N)}]$ is bounded by a constant independent of $N$, we may employ the Ito isometry to find a constant $C_2>0$ depending only on $T$ and $\sup_{u\in[0,1]}\sigma(u)$ such that for all $N\in\nn$ and $A>0$:
\begin{eqnarray*}
Q^{(N)}\Big(\sup_{t\in[0,T]} (\rho^{(N)}(t),\phi)\geq A+C_1 T\Big)\leq\frac{C_2}{A}.
\end{eqnarray*}
Hence, we can choose $A$ such that $\frac{C_2}{A}<\epsilon$ and let 
\begin{eqnarray}
K_0=\{\alpha\in M_1(\rr)|\;(\alpha,\phi)\leq A+C_1 T\}.
\end{eqnarray}
As explained in the proof of Lemma 1.4 in \cite{ga} the compactness of the set $K_0$ in $M_1(\rr)$ is a consequence of Prokhorov's Theorem. It remains to show that $\ev^{Q^{(N)}}[b^{(N)}]$ is bounded by a constant independent of $N$. It is clear from the definition of $\rho^{(N)}(0)$, $N\in\nn$ that $\ev^{Q^{(N)}}[b^{(N)}]$ is finite for all $N\in\nn$. Moreover, for $N\geq4$ we can use $\phi(x)\leq 1+|x|$ and inequality (\ref{plb}) to compute
\begin{eqnarray*}
&&\ev^{Q^{(N)}}[b^{(N)}]=\frac{1}{N}\sum_{i=1}^N \ev[\phi(Y_i(0)-Y_{M(N)}(0))]\leq 1+\frac{1}{N}\sum_{i=1}^N \ev[|Y_i(0)-Y_{M(N)}(0)|]\\
&&=1+\frac{1}{N}\sum_{i=1}^{M(N)-1}\sum_{j=i}^{M(N)-1}\frac{1}{a_j^{(N)}}+\frac{1}{N}\sum_{i=M(N)+1}^N\sum_{j=M(N)}^{i-1}\frac{1}{a_j^{(N)}}\\
&&\leq1+\frac{2}{N}\cdot\frac{\sup_{u\in[0,1]}\sigma(u)^2}{\omega_0}\sum_{i=1}^{M(N)-1}\sum_{j=i}^{M(N)-1}\frac{N}{j(N-j)}.
\end{eqnarray*}
Finally, the upper bound
\begin{eqnarray*}
\sum_{i=1}^{M(N)-1}\sum_{j=i}^{M(N)-1}\frac{1}{j(N-j)}=\sum_{j=1}^{M(N)-1} \frac{1}{N-j}\leq\log(N-1)-\log(N-M(N))\leq\log 3
\end{eqnarray*}
for all $N\geq4$ shows that $\ev^{Q^{(N)}}[b^{(N)}]$ is bounded by a constant independent of $N$.\\\\
2) To prove the existence of sets $K_1,K_2,\dots$ with the desired properties it suffices to show that for any fixed $r\in\nn$ the sequence of probability measures $P^{(N),f_r}$, $N\in\nn$ on $C([0,T])$ induced by ${\mathcal R}^{(N)}$, $N\in\nn$ through the mapping $\xi\mapsto(\xi(.),f_r)$ is tight. To this end, we fix an $r\in\nn$ and aim to deduce the tightness of the sequence $P^{(N),f_r}$, $N\in\nn$ from Theorem 1.3.2 of \cite{sv}. To do this we need to show 
\begin{eqnarray*}
\lim_{\theta\uparrow\infty}\inf_{N\in\nn} P^{(N),f_r}(|y(0)|\leq\theta)=1
\end{eqnarray*} 
and
\begin{eqnarray*}
\forall\Delta>0:\quad 
\lim_{\epsilon\downarrow0}\limsup_{N\rightarrow\infty} P^{(N),f_r}\Big(\sup_{0\leq s\leq t\leq T,t-s\leq\epsilon}|y(t)-y(s)|>\Delta\Big)=0.
\end{eqnarray*}
The first assertion follows immediately by considering $\theta>\sup_{x\in\rr}|f_r(x)|$. To show the second assertion we fix a $\Delta>0$, define $Z(t)$, $t\in[0,T]$ as in step 1, but replacing $\phi$ by $f_r$ and redefine the constant $C_1$ correspondingly. Using the $L^2$-version of Doob's maximal inequality for non-negative continuous submartingales we obtain for each $0<\epsilon<\frac{\Delta}{C_1}$:
\begin{eqnarray*}
P^{(N),f_r}\Big(\sup_{0\leq s\leq t\leq T,t-s\leq\epsilon}|y(t)-y(s)|>\Delta\Big)
\leq Q^{(N)}\Big(\sup_{0\leq s\leq t\leq T,t-s\leq\epsilon}|Z(t)-Z(s)|>\Delta-C_1\epsilon\Big)\\
\leq \sum_{l=0}^{\left\lfloor \frac{T}{\epsilon}\right\rfloor-1}
Q^{(N)}\Big(\sup_{l\epsilon\leq s\leq \min((l+2)\epsilon,T)}|Z(s)-Z(l\epsilon)|\geq\frac{\Delta-C_1\epsilon}{2}\Big)\\
\leq \sum_{l=0}^{\left\lfloor\frac{T}{\epsilon}\right\rfloor-1}
\left(\frac{\Delta-C_1\epsilon}{2}\right)^{-2}\ev^{Q^{(N)}}\Big[|Z(\min((l+2)\epsilon,T))-Z(l\epsilon)|^2\Big]
\end{eqnarray*} 
where $\left\lfloor.\right\rfloor$ denotes the integer part of a positive real number. The Ito isometry shows that the latter expression is bounded by $\frac{T}{\epsilon}\left(\frac{\Delta-C_1\epsilon}{2}\right)^{-2}\frac{C_3\epsilon}{N}$ with a constant $C_3>0$ depending only on $\sup_{x\in\rr}|f_r'(x)|$ and $\sup_{u\in[0,1]}\sigma(u)$. Taking limits we end up with the second assertion. We conclude that the sequence ${\mathcal R}^{(N)}$, $N\in\nn$ is tight.\\\\
3) To prove (\ref{mcv2}) we let ${\mathcal R}^\infty$ be the limit of a converging subsequence ${\mathcal R}^{(N_k)}$, $k\in\nn$ of the sequence ${\mathcal R}^{(N)}$, $N\in\nn$. Next, for each $k\in\nn$ we let $\widetilde{\rho}^{(N_k)}$ be a random variable with distribution ${\mathcal R}^{(N_k)}$ and $\widetilde{\rho}^\infty$ be a random variable with distribution ${\mathcal R}^\infty$, all defined on the same probability space and such that $\widetilde{\rho}^{(N_k)}\rightarrow_{k\rightarrow\infty}\widetilde{\rho}^\infty$ in $C([0,T],M_1(\rr))$ almost surely. This is possible due to the Skorohod Representation Theorem in the form of Theorem 3.5.1 in \cite{du}. Indeed, the metric space $C([0,T],M_1(\rr))$ is separable, since the countable set of functions whose values at $0,\frac{1}{j},\frac{2}{j},\dots,1$ belong to a fixed countable dense subset of $M_1(\rr)$ and which interpolate linearly on the intervals  $\Big[0,\frac{1}{j}\Big],\dots,\Big[\frac{j-1}{j},1\Big]$ for a $j\in\nn$ is dense in $C([0,T],M_1(\rr))$.\\\quad\\ 
From the dynamics computed in step 1 we observe that for each $k\in\nn$ and $f\in C_c^3(\rr)$ the process $(\widetilde{\rho}^{(N_k)}(t),f)$, $t\in[0,T]$ is a semimartingale in the sense of definition II.2.6 in \cite{js} with its characteristics being given by
\begin{eqnarray*}
&&\int_0^t \int_\rr f'(x)\mu(F_{\widetilde{\rho}^{(N_k)}(s)}(x))+\frac{1}{2}f''(x)\sigma(F_{\widetilde{\rho}^{(N_k)}(s)}(x))^2\;\widetilde{\rho}^{(N_k)}(s)(dx)\;ds,\\
&&\frac{1}{N_k}\int_0^t\int_\rr f'(x)^2 \sigma(F_{\widetilde{\rho}^{(N_k)}(s)}(x))^2\;\widetilde{\rho}^{(N_k)}(s)(dx)\;ds
\end{eqnarray*}
which we denote by $B^{(N_k)}(t)$ and $C^{(N_k)}(t)$, respectively. We claim that in order to establish (\ref{mcv2}) it suffices to show that for any $f\in C_c^3(\rr)$ we have
\begin{eqnarray}
&&\sup_{t\in[0,T]}\Big|B^{(N_k)}(t)-\int_0^t(\widetilde{\rho}^\infty(s),L_{\widetilde{\rho}^\infty(s)}f)\;ds\Big|\stackrel{p}{\rightarrow}0, \label{conv1}\\
&&\ev[C^{(N_k)}(T)]\rightarrow_{k\rightarrow\infty}0 \label{conv2}
\end{eqnarray}
where $\stackrel{p}{\rightarrow}$ denotes the convergence in probability. Indeed, the second convergence together with the $L^2$-version of Doob's maximal inequality for non-negative continuous submartingales would imply
\begin{eqnarray}
\sup_{t\in[0,T]}\Big|(\widetilde{\rho}^{(N_k)}(t),f)-(\widetilde{\rho}^{(N_k)}(0),f)-B^{(N_k)}(t)\Big|\stackrel{p}{\rightarrow}0.
\end{eqnarray}
Hence, from the first convergence we would be able to conclude that
\begin{eqnarray}
\sup_{t\in[0,T]}\Big|(\widetilde{\rho}^{(N_k)}(t),f)-(\widetilde{\rho}^{(N_k)}(0),f)-\int_0^t(\widetilde{\rho}^\infty(s),L_{\widetilde{\rho}^\infty(s)}f\;ds\Big|\stackrel{p}{\rightarrow}0.
\end{eqnarray}
By a diagonalization argument relying on the separability of $C_c^3(\rr)$, endowed with the topology of uniform convergence of functions and their first and second derivatives, we would be able to find a subsequence of $N_k$, $k\in\nn$ such that the latter convergence holds for all $f\in C_c^3(\rr)$ in the almost sure sense. But since  $(\widetilde{\rho}^{(N_k)}(t),f)-(\widetilde{\rho}^{(N_k)}(0),f)-
\int_0^t(\widetilde{\rho}^\infty(s),L_{\widetilde{\rho}^\infty(s)}f)\;ds$ converges to $(\widetilde{\rho}^\infty(t),f)-(\widetilde{\rho}^\infty(0),f)-\int_0^t(\widetilde{\rho}^\infty(s),L_{\widetilde{\rho}^\infty(s)}f)\;ds$
in the limit $k\rightarrow\infty$ for all $f\in C_c^3(\rr)$ and $t\in[0,T]$ almost surely, we would obtain equation (\ref{mcv2}).\\\\
4) We now show the two claimed convergence results. The convergence in (\ref{conv2}) is a direct consequence of the boundedness of $f'$ and $\sigma$. To prove the convergence in (\ref{conv1}) we introduce for any $\epsilon>0$ and $x\in\rr$ a Lipschitz function $f_x^\epsilon$ such that 
\begin{eqnarray}
1_{(-\infty,x-\epsilon]}\leq f^\epsilon_x\leq 1_{(-\infty,x]}
\end{eqnarray}
and $f^\epsilon_{x'}$ is a translate of $f^\epsilon_x$ by $x'-x$ for any $x,x'\in\rr$ and $\epsilon>0$. Moreover, for any probability measure $\alpha$ on $\rr$ we set $F^\epsilon_\alpha(x)=(\alpha,f^\epsilon_x)$. We note $B^{(N_k)}(t)=\int_0^t(\widetilde{\rho}^{(N_k)}(s),L_{\widetilde{\rho}^{(N_k)}(s)}f)\;ds$ for all $t\in[0,T]$ and $k\in\nn$ and deduce from the triangle inequality that
\begin{eqnarray*}
&&\sup_{t\in[0,T]}\Big|\int_0^t(\widetilde{\rho}^{(N_k)}(s),L_{\widetilde{\rho}^{(N_k)}(s)}f)\;ds-\int_0^t(\widetilde{\rho}^\infty(s),L_{\widetilde{\rho}^\infty(s)}f)\;ds\Big|\\
&&\leq\sup_{t\in[0,T]}\Big|\int_0^t(\widetilde{\rho}^{(N_k)}(s),f'(\mu\circ F_{\widetilde{\rho}^{(N_k)}(s)}))
-(\widetilde{\rho}^\infty(s),f'(\mu\circ F_{\widetilde{\rho}^\infty(s)}))\;ds\Big|\\
&&+\sup_{t\in[0,T]}\Big|\int_0^t(\widetilde{\rho}^{(N_k)}(s),\frac{f''}{2}(\sigma\circ F_{\widetilde{\rho}^{(N_k)}(s)})^2)
-(\widetilde{\rho}^\infty(s),\frac{f''}{2}(\sigma\circ F_{\widetilde{\rho}^\infty(s)})^2)\;ds\Big|
\end{eqnarray*}
where $\circ$ denotes the composition of functions. We claim that the latter two terms converge to zero in probability. Since the proof of this claim is identical for both terms, we only carry it out for the first (drift) term. To this end, we fix an $\epsilon>0$ and a $k\in\nn$ and observe  
\begin{eqnarray*}
&&\sup_{t\in[0,T]}\Big|\int_0^t(\widetilde{\rho}^{(N_k)}(s),f'(\mu\circ F_{\widetilde{\rho}^{(N_k)}(s)}))
-(\widetilde{\rho}^\infty(s),f'(\mu\circ F_{\widetilde{\rho}^\infty(s)}))\;ds\Big|\\
&&\leq\int_0^T |(\widetilde{\rho}^{(N_k)}(s),f'(\mu\circ F_{\widetilde{\rho}^{(N_k)}(s)}))-(\widetilde{\rho}^{(N_k)}(s),f'(\mu\circ F^\epsilon_{\widetilde{\rho}^{(N_k)}(s)}))|\;ds\\
&&+\sup_{t\in[0,T]}\Big|\int_0^t(\widetilde{\rho}^{(N_k)}(s),f'(\mu\circ F^\epsilon_{\widetilde{\rho}^{(N_k)}(s)}))
-(\widetilde{\rho}^\infty(s),f'(\mu\circ F^\epsilon_{\widetilde{\rho}^\infty(s)}))\;ds\Big|\\
&&+\int_0^T |(\widetilde{\rho}^\infty(s),f'(\mu\circ F^\epsilon_{\widetilde{\rho}^\infty(s)}))-(\widetilde{\rho}^\infty(s),f'(\mu\circ F_{\widetilde{\rho}^\infty(s)}))|\;ds.
\end{eqnarray*}
We call the summands on the right-hand side (I), (II) and (III). We will bound the three terms consecutively.\\\quad\\ 
Denoting by $Y^R_1(t)\leq\dots\leq Y^R_{N_k}(t)$ the ordered particles of the normalized system $R_1(t),\dots,R_{N_k}(t)$ for any time $t\in[0,T]$ we can bound term (I) from above by
\begin{eqnarray*}
C_4\int_0^T \frac{1}{N_k}\sum_{i=1}^{N_k}\frac{1}{N_k}\Big|\{1\leq j\leq i|\;Y^R_j(s)\in(Y^R_i(s)-\epsilon,Y^R_i(s)]\}\Big|\;ds
\end{eqnarray*}
where the constant $C_4$ is the product of $\sup_{x\in\rr}|f'(x)|$ and $\sup_{u\in[0,1]}|\mu'(u)|$. By Fubini's Theorem and the defining property of the initial condition (see the paragraph preceeding Theorem 1.1) it follows that the expectation of (I) is bounded above by $\frac{C_4T}{N_k^2}\sum_{i=1}^{N_k}\sum_{j=1}^i\pp(Y_i(0)-Y_j(0)<\epsilon)$. To bound this expression further we first choose a $C_5>0$ such that $a_i^{(N)}\leq C_5 N$ for all $N\in\nn$ and $1\leq i\leq N-1$. This is possible due to the obvious bound $a_i^{(N)}\leq\frac{4i(N-i)}{N}\cdot\frac{\sup_{u\in[0,1]}|\mu(u)|}{\inf_{u\in[0,1]}\sigma(u)^2}$ which holds for all $N\in\nn$ and $1\leq i\leq N-1$. Next, for any fixed $k\in\nn$ we let $E_1,\dots,E_{N_k-1}$ be i.i.d. exponential random variables with parameter $C_5 N_k$ and $P$ be a Poisson random variable with parameter $C_5 N_k\epsilon$. Then from the scaling property of exponential random variables we deduce
\begin{eqnarray*}
&&\frac{C_4T}{N_k^2}\sum_{i=1}^{N_k}\sum_{j=1}^i\pp(Y_i(0)-Y_j(0)\leq\epsilon)
\leq\frac{C_4T}{N_k^2}\sum_{i=1}^{N_k}\sum_{j=1}^i\pp(E_j+\dots+E_{i-1}\leq\epsilon)\\
&&=\frac{C_4T}{N_k^2}\sum_{i=1}^{N_k}\sum_{j=1}^i\pp(P\geq i-j)
\leq \frac{C_4T}{N_k}\Big(1+\ev[P]\Big)=\frac{C_4T}{N_k}\Big(1+C_5 N_k\epsilon\Big). 
\end{eqnarray*}   
All in all, we conclude that
\begin{eqnarray*}
\limsup_{k\rightarrow\infty}\ev\left[\int_0^T |(\widetilde{\rho}^{(N_k)}(s),f'(\mu\circ F_{\widetilde{\rho}^{(N_k)}(s)}))
-(\widetilde{\rho}^{(N_k)}(s),f'(\mu\circ F^\epsilon_{\widetilde{\rho}^{(N_k)}(s)}))|\;ds\right]
\end{eqnarray*}
is bounded above by $C_4C_5T\epsilon$.\\\quad\\ 
Term (II) can be estimated from above by
\begin{eqnarray*}
&&\sup_{t\in[0,T]}\Big|\int_0^t(\widetilde{\rho}^{(N_k)}(s),f'(\mu\circ F^\epsilon_{\widetilde{\rho}^{(N_k)}(s)}))
-(\widetilde{\rho}^{(N_k)}(s),f'(\mu\circ F^\epsilon_{\widetilde{\rho}^\infty(s)}))\;ds\Big|\\
&&+\sup_{t\in[0,T]}\Big|\int_0^t(\widetilde{\rho}^{(N_k)}(s),f'(\mu\circ F^\epsilon_{\widetilde{\rho}^\infty(s)}))
-(\widetilde{\rho}^\infty(s),f'(\mu\circ F^\epsilon_{\widetilde{\rho}^\infty(s)}))\;ds\Big|.
\end{eqnarray*}
Now, recalling the definition of $C_4$ we bound the first summand from above by
\begin{eqnarray} \label{c4b}
C_4\sup_{t\in[0,T]}\int_0^t \sup_{x\in\rr}|F^\epsilon_{\widetilde{\rho}^{(N_k)}(s)}(x)-F^\epsilon_{\widetilde{\rho}^\infty(s)}(x)|\;ds. 
\end{eqnarray}
Since for each $x\in\rr$ the function $f^\epsilon_x$ is Lipschitz with the Lipschitz constant being independent of $x$ and the convergence $d(\widetilde{\rho}^{(N_k)}(s),\widetilde{\rho}^\infty(s))\rightarrow_{k\rightarrow\infty}0$ is uniform in $s$, the expression in (\ref{c4b}) tends to zero almost surely in the limit $k\rightarrow\infty$. Moreover, for any $x,x'\in\rr$ it holds
\begin{eqnarray*}
|F^\epsilon_{\widetilde{\rho}^\infty(s)}(x)-F^\epsilon_{\widetilde{\rho}^\infty(s)}(x')|\leq\sup_{y\in\rr} |f^\epsilon_x(y)-f^\epsilon_x(y-x+x')|
\leq Lip(f^\epsilon_0)\cdot|x-x'|. 
\end{eqnarray*}
Hence, the uniformity in $s$ of the convergence $d(\widetilde{\rho}^{(N_k)}(s),\widetilde{\rho}^\infty(s))\rightarrow_{k\rightarrow\infty}0$ implies that the second summand in the bound on term (II) tends to zero almost surely in the limit $k\rightarrow\infty$. Thus, from the Dominated Convergence Theorem we deduce that the expectation of term (II) converges to zero in the limit $k\rightarrow\infty$. \\\quad\\ 
Due to the inequality
\begin{eqnarray}
|F^\epsilon_{\widetilde{\rho}^\infty(s)}(x)-F_{\widetilde{\rho}^\infty(s)}(x)|\leq F^\epsilon_{\widetilde{\rho}^\infty(s)}(x+\epsilon)-F^\epsilon_{\widetilde{\rho}^\infty(s)}(x)
\end{eqnarray}
for all $x\in\rr$ the expectation of term (III) can be bounded above by
\begin{eqnarray*}
C_4\;\ev\left[\int_0^T \int_\rr (F^\epsilon_{\widetilde{\rho}^\infty(s)}(x+\epsilon)-F^\epsilon_{\widetilde{\rho}^\infty(s)}(x)) \widetilde{\rho}^\infty(s)(dx)\;ds\right].
\end{eqnarray*}
Moreover, applying the triangle inequality as in the proof of the upper bound on term (II) one shows that with probability $1$ the integrand in the $ds$-integral is the limit of the corresponding objects with $\widetilde{\rho}^\infty$ replaced by $\widetilde{\rho}^{(N_k)}$. Indeed, the difference between the two is bounded above in absolute value by $4(Lip(f^\epsilon_0)+1)\cdot d(\widetilde{\rho}^{(N_k)}(s),\widetilde{\rho}^\infty(s))$. We conclude from Fatou's Lemma that the expectation of term (III) is bounded above by  
\begin{eqnarray*}
C_4\liminf_{k\rightarrow\infty}\ev\left[\int_0^T \int_\rr (F^\epsilon_{\widetilde{\rho}^{(N_k)}(s)}(x+\epsilon)-F^\epsilon_{\widetilde{\rho}^{(N_k)}(s)}(x)) \widetilde{\rho}^{(N_k)}(s)(dx)\;ds\right].
\end{eqnarray*}
Bounding the integrand in the $ds$-integral in the same way as in the corresponding estimate on term (I) one obtains the upper bound
\begin{eqnarray*}
C_4\liminf_{k\rightarrow\infty}
\ev\left[\int_0^T \frac{1}{N_k}\sum_{i=1}^{N_k}\frac{1}{N_k}\Big|\{1\leq j\leq N_k|\;\;|Y^R_j(s)-Y^R_i(s)|\leq\epsilon\}\Big|\;ds\right].
\end{eqnarray*}
Finally, proceeding as in the upper bound on term (I) we deduce that the expectation of term (III) is bounded above by $2C_4C_5T\epsilon$.\\\quad\\
All in all, we have shown that
\begin{eqnarray*}
\limsup_{k\rightarrow\infty}\ev\left[\sup_{t\in[0,T]}\Big|\int_0^t(\widetilde{\rho}^{(N_k)}(s),f'(\mu\circ F_{\widetilde{\rho}^{(N_k)}(s)}))
-(\widetilde{\rho}^\infty(s),f'(\mu\circ F_{\widetilde{\rho}^\infty(s)}))ds\Big|\right]
\end{eqnarray*}
is bounded above by $C_6\epsilon$ with $C_6=3C_4C_5T$. By taking the limit $\epsilon\downarrow0$ we deduce that the term inside the latter expectation tends to zero in $L^1$ in the limit $k\rightarrow\infty$ and so, in particular, it converges to zero in probability. \ep 

\section{Proof of Theorem 1.2}

To be ready to prove Theorem 1.2 we show next that the initial probability measure $\xi(0)$ is the same under each accumulation point ${\mathcal R}^\infty$ of the set $\Xi$ provided that the functions $\mu$ and $\sigma^2$ are affine.
{\prop Let the functions $\mu$ and $\sigma^2$ be affine, so that $\sigma^2(u)=cu+d$ for some constants $c,d\in\rr$. Then for any choice of ${\mathcal R}^\infty$ as in Theorem 1.1 the distribution of the initial probability measure $\xi(0)$ under ${\mathcal R}^\infty$ is a Dirac probability measure. Moreover, the atom of the latter is given by the unique $\lambda^\infty\in M_1(\rr)$ whose quantiles $q^\infty(u)$, $u\in(0,1)$ are given by 
\begin{eqnarray}
q^\infty(u)=-\frac{c+d}{|\mu'(0)|}\log(2-2u)+\frac{d}{|\mu'(0)|}\log(2u).
\end{eqnarray}
In particular, $\lambda^\infty$ is absolutely continuous with respect to the Lebesgue measure on $\rr$. \\\quad\\}
{\it Proof.} 1) For each $N\in\nn$ let $q^{(N)}(u)$, $u\in(0,1)$ be the (random) quantiles of the initial probability measure $\xi(0)$ under ${\mathcal R}^{(N)}$. Then for any fixed $\frac{1}{2}<u<1$ it holds
\begin{eqnarray*}
&&q^{(N)}(u)=\sum_{i=M(N)}^{\left\lceil uN\right\rceil-1} (Y_{i+1}(0)-Y_i(0))\\
&&=\frac{\sum_{i=M(N)}^{\left\lceil uN\right\rceil-1} (Y_{i+1}(0)-Y_i(0))-m^{(N)}}{\sqrt{v^{(N)}}}\cdot\sqrt{v^{(N)}}+m^{(N)}
\end{eqnarray*}
for all $N\in\nn$ such that $\left\lceil uN\right\rceil-1\geq M(N)$ where 
\begin{eqnarray}
m^{(N)}=\sum_{i=M(N)}^{\left\lceil uN\right\rceil-1}\frac{1}{a_i^{(N)}},\quad 
v^{(N)}=\sum_{i=M(N)}^{\left\lceil uN\right\rceil-1} \frac{1}{(a_i^{(N)})^2}
\end{eqnarray}
and $\left\lceil x\right\rceil$ denotes the smallest integer greater or equal to $x$ for any $x\in\rr$. Using the definition of $a_i^{(N)}$, $1\leq i\leq N-1$, $N\in\nn$ and the assumption that $\mu$ is affine we see that for any $\frac{1}{2}<u<1$ there exist constants $\omega_1(u),\omega_2(u)>0$ such that $\omega_1(u)N\leq a_i^{(N)}\leq \omega_2(u)N$ for all $M(N)\leq i\leq \left\lceil uN\right\rceil-1$, $N\in\nn$. From this and Lindeberg's Central Limit Theorem we deduce that for any fixed $\frac{1}{2}<u<1$ it holds
\begin{eqnarray}
\lim_{N\rightarrow\infty} q^{(N)}(u)=\lim_{N\rightarrow\infty} m^{(N)}
\end{eqnarray}
in distribution, provided that the latter limit exists. Indeed, we can check Lindeberg's condition for the triangular array $\frac{Y_{i+1}(0)-Y_i(0)-\frac{1}{a_i^N}}{\sqrt{v^{(N)}}}$, $M(N)\leq i\leq \left\lceil uN\right\rceil-1$, $N\in\nn$ as follows: for any $\epsilon>0$ and any $N\in\nn$ such that $\left\lceil uN\right\rceil-1\geq M(N)$ we have the estimates
\begin{eqnarray*}
&&\frac{1}{v^{(N)}}\cdot\sum_{i=M(N)}^{\left\lceil uN\right\rceil-1}\ev\Big[\Big(Y_{i+1}(0)-Y_i(0)-\frac{1}{a_i^{(N)}}\Big)^2\cdot 1_{|Y_{i+1}(0)-Y_i(0)-\frac{1}{a_i^{(N)}}|\geq\epsilon\sqrt{v^{(N)}}}\Big]\\
&&\leq \frac{1}{v^{(N)}}\cdot\sum_{i=M(N)}^{\left\lceil uN\right\rceil-1}\frac{1}{(a_i^{(N)})^2}\cdot
\ev\Big[G_{i,N}^2\cdot 1_{|G_{i,N}|\geq\epsilon\sqrt{v^{(N)}}a_i^{(N)}}\Big]\\
&&\leq \ev\Big[G_{M(N),N}^2\cdot 1_{|G_{M(N),N}|\geq\epsilon\omega_3(u)\sqrt{N}}\Big]
\end{eqnarray*}
where we have set $G_{i,N}=(Y_{i+1}(0)-Y_i(0))a_i^{(N)}-1$ for $M(N)\leq i\leq \left\lceil uN\right\rceil-1$, $N\in\nn$ and $\omega_3(u)>0$ is a constant such that $\sqrt{v^{(N)}}a_i^{(N)}\geq \omega_3(u)\sqrt{N}$ for all $M(N)\leq i\leq \left\lceil uN\right\rceil-1$ and $N\in\nn$. The latter expectation tends to zero in the limit $N\rightarrow\infty$, since the random variable $G_{M(N),N}+1$ is distributed according to the exponential distribution with parameter $1$ for all $N\geq3$. Hence, Lindeberg's condition is satisfied.\\\quad\\ 
2) Now, plugging in the definition of $a_i^{(N)}$ for $M(N)\leq i\leq \left\lceil uN\right\rceil-1$ and $N\in\nn$ we calculate
\begin{eqnarray*}
&&\lim_{N\rightarrow\infty} m^{(N)}
=\lim_{N\rightarrow\infty}\sum_{i=M(N)}^{\left\lceil uN\right\rceil-1}\frac{N}{i(N-i)}\cdot\frac{\frac{ci+c/2}{N}+d}{|\mu'(0)|}\\
&&=\frac{c}{|\mu'(0)|}\cdot\lim_{N\rightarrow\infty}\sum_{i=M(N)}^{\left\lceil uN\right\rceil-1}\frac{1}{N-i}
+\frac{d}{|\mu'(0)|}\cdot\lim_{N\rightarrow\infty}\sum_{i=M(N)}^{\left\lceil uN\right\rceil-1}\left(\frac{1}{i}+\frac{1}{N-i}\right)\\
&&=\lim_{N\rightarrow\infty}\Big(\frac{c+d}{|\mu'(0)|}\log\Big(\frac{N-M(N)}{N-uN}\Big)+\frac{d}{|\mu'(0)|}\log\Big(\frac{uN}{M(N)}\Big)\Big).
\end{eqnarray*}
Recalling that $\frac{N}{2}\leq M(N)\leq\frac{N}{2}+1$ we can compute the last limit to 
\begin{eqnarray*}
-\frac{c+d}{|\mu'(0)|}\log(2-2u)+\frac{d}{|\mu'(0)|}\log(2u).
\end{eqnarray*}
An analogous application of Lindeberg's Central Limit Theorem and similar calculations to the ones above show for the case $0<u\leq\frac{1}{2}$:
\begin{eqnarray*}
\lim_{N\rightarrow\infty} q^{N}(u)=-\lim_{N\rightarrow\infty} \sum_{i=\left\lceil uN\right\rceil}^{M(N)-1} \frac{1}{a_i^{(N)}}
=-\frac{c+d}{|\mu'(0)|}\log(2-2u)+\frac{d}{|\mu'(0)|}\log(2u)
\end{eqnarray*}
in distribution.\\\quad\\
3) If ${\mathcal R}^\infty$ is as in Theorem 1.1, then the Skorohod Representation Theorem in the form of Theorem 3.5.1 in \cite{du} shows that we can find an increasing sequence $N_k$, $k\in\nn$ of natural numbers and random variables $\widetilde{\rho}^{(N_k)}(0)$, $k\in\nn$ and $\widetilde{\rho}^\infty(0)$ defined on the same probability space such that for each $k\in\nn$ the distribution of the random variable $\widetilde{\rho}^{(N_k)}(0)$ is given by the law of $\xi(0)$ under ${\mathcal R}^{(N_k)}$, $\widetilde{\rho}^\infty(0)$ is distributed according to the law of $\xi(0)$ under ${\mathcal R}^\infty$ and $\widetilde{\rho}^{(N_k)}(0)\rightarrow_{k\rightarrow\infty}\widetilde{\rho}^\infty(0)$ weakly with probability $1$. It follows that the quantile functions of $\widetilde{\rho}^{(N_k)}(0)$ converge in the limit $k\rightarrow\infty$ to the quantile function of $\widetilde{\rho}^\infty(0)$ at all continuity points of the latter almost surely (see e.g. the proof of Theorem 2.2.2 in \cite{dur}). By Fubini's Theorem we obtain that the $u$-quantile of $\widetilde{\rho}^{(N_k)}(0)$ converges to the $u$-quantile of $\widetilde{\rho}^\infty(0)$ almost surely in the limit $k\rightarrow\infty$ for Lebesgue almost every $u\in(0,1)$ and, in particular, for all $u$ in a countable dense subset of $(0,1)$. Due to the monotonicity of quantile functions and the computations in steps 1 and 2 the quantile function of $\widetilde{\rho}^\infty(0)$ has to coincide with $q^\infty$ (defined in the statement of the proposition) with probability $1$. Hence, the distribution of $\widetilde{\rho}^\infty(0)$, which is the same as the law of $\xi(0)$ under ${\mathcal R}^\infty$, is given by the Dirac probability measure described in the proposition. Finally, the probability measure $\lambda^\infty$ is absolutely continuous with respect to the Lebesgue measure on $\rr$, because its quantile function $q^\infty$ is continuously differentiable and strictly increasing on $(0,1)$. \ep\\\quad\\ 
Combining the ideas of \cite{ga} and \cite{jo} with a result in \cite{gi} we can now prove Theorem 1.2. In the proof we use the following notations. For a measurable subset $S$ of a Euclidean space we write $L^p(S)$ and $\|.\|_{L^p(S)}$ for the space of functions $f:S\rightarrow\rr$ such that $|f|^p$ is integrable with respect to the restriction of the Lebesgue measure to $S$ and the corresponding $L^p$-norm, respectively, where $p$ is a real number in $[1,\infty)$. In addition, for any real-valued random variable $Y$ we denote by ${\mathcal L}(Y)$ the law of $Y$. \\\quad\\ 
{\it Proof of Theorem 1.2.} 1) It suffices to prove that the initial value problem (\ref{mcv}), (\ref{mcvic}) has at most one solution in $C([0,T],M_1(\rr))$ for each initial condition which is absolutely continuous with respect to the Lebesgue measure on $\rr$, since then the law of large numbers follows from Theorem 1.1 and Proposition 3.1. To achieve the former we fix an absolutely continuous initial condition $\lambda$, let $\nu_i$, $i\in\{1,2\}$ be two solutions of (\ref{mcv}), (\ref{mcvic}) in $C([0,T],M_1(\rr))$ with this initial condition and will prove $\nu_1=\nu_2$. We show first that for $i\in\{1,2\}$ the measures $\nu_i(t)$, $t\in[0,T]$ are given by one-dimensional distributions of solutions of appropriate martingale problems. To this end, for $i\in\{1,2\}$ we define $Z_i(t)$, $t\in[0,T]$ as the respective unique (in law) solutions of the martingale problems associated with the families of operators $L_{\nu_i(t)}$, $t\in[0,T]$ such that for any $f\in C^\infty_c(\rr)$ the processes
\begin{eqnarray*}
f(Z_i(t))-f(Z_i(0))-\int_0^t (L_{\nu_i(s)}f)(Z_i(s))\;ds,\quad t\in[0,T],
\end{eqnarray*}
$i\in\{1,2\}$ are martingales and ${\mathcal L}(Z_1(0))={\mathcal L}(Z_2(0))=\lambda$. Due to Exercise 7.3.3 in \cite{sv} the processes $Z_1(t)$, $t\in[0,T]$ and $Z_2(t)$, $t\in[0,T]$ are well-defined and ${\mathcal L}(Z_i(.))\in C([0,T],M_1(\rr))$ for $i\in\{1,2\}$. We claim that $\nu_i(.)={\mathcal L}(Z_i(.))$ for $i\in\{1,2\}$. To prove the claim we fix an $i\in\{1,2\}$. By their respective definitions $\nu_i(t)$, $t\in[0,T]$ and ${\mathcal L}(Z_i(t))$, $t\in[0,T]$ solve the initial value problem for the Fokker-Planck equation 
\begin{eqnarray}
&&\forall f\in C_c^\infty(\rr):\quad (\xi(t),f)-(\xi(0),f)=\int_0^t (\xi(s),L_{\nu_i(s)}f)\;ds, \label{fp}\\
&&\xi(0)=\lambda
\end{eqnarray}
on $[0,T]$. We show now that the latter has a unique solution in $C([0,T],M_1(\rr))$. To this end, we observe that the operator ${\mathcal R}=\frac{\partial}{\partial t}+L_{\nu_i(.)}$ is continuous as an operator from the Sobolev space $W^{1,2,p}([0,T]\times[-r,r])$ (the space functions in $L^p([0,T]\times[-r,r])$ whose generalized first time derivative and generalized first two spatial derivatives belong to $L^p([0,T]\times[-r,r])$, endowed with the usual Sobolev norm) into $L^p([0,T]\times[-r,r])$ for any $p\geq1$ and any $r>0$. This is due to the boundedness of $\mu$ and $\sigma$. Moreover, following the steps in the proof of Theorem A.1 in \cite{ga} we obtain for any $p>3$ and $r>0$:
\begin{eqnarray}
\quad\quad\;\;\int_0^T \Big(\nu_i(t),\Big(\frac{\partial}{\partial t}+L_{\nu_i(t)}\Big)f\Big)\;dt
\geq \int_0^T \Big({\mathcal L}(Z_i(t)),\Big(\frac{\partial}{\partial t}+L_{\nu_i(t)}\Big)f\Big)\;dt
\end{eqnarray} 
for all $f\in W_0^{1,2,p}([0,T]\times[-r,r])$ (the space of functions in $W^{1,2,p}([0,T]\times[-r,r])$ vanishing on $([0,T]\times\{-r,r\})\cup(\{T\}\times[-r,r])$) such that ${\mathcal R}f\geq0$ Lebesgue almost everywhere. Hereby, we have used the convention $(\frac{\partial}{\partial t}+L_{\nu_i(t)})f=0$ on the complement of $[0,T]\times[-r,r]$ in $[0,T]\times\rr$. By Theorem 9.1 in chapter IV of \cite{ls} the image of the just described functions under ${\mathcal R}$ is given by $L^p_+([0,T]\times[-r,r])$, the set of functions in $L^p([0,T]\times[-r,r])$ which are non-negative Lebesgue almost everywhere. Indeed, from Proposition 3.1 and the assumption of the theorem we conclude that for every $t\in[0,T]$ the measure $\nu_i(t)$ is absolutely continuous with respect to the Lebesgue measure on $\rr$. Thus, the function $(t,x)\mapsto\sigma(F_{\nu_i(t)}(x))$ is continuous and Theorem 9.1 in chapter IV of \cite{ls} is applicable. Since $r>0$ was arbitrary, we deduce from the Monotone Convergence Theorem that
\begin{eqnarray*} 
\int_0^T \int_\rr g(x)\cdot 1_{[0,t]}(s)\;\nu_i(s)(dx)\;ds
\geq \int_0^T \int_\rr g(x)\cdot 1_{[0,t]}(s)\;{\mathcal L}(Z_i(s))(dx)\;ds
\end{eqnarray*}
for all non-negative continuous bounded functions $g$ on $\rr$ and all $t\in[0,T]$. Using the same inequality with $(\sup_{x\in\rr}g(x))-g$ instead of $g$ we infer that equality must hold in the latter inequality. Differentiating the resulting identity with respect to $t$ we see that $\nu_i(t)={\mathcal L}(Z_i(t))$ for all $t\in[0,T]$. \\\quad\\
2) To finish the proof we aim to apply Theorem 4 of \cite{gi}. Approximating the functions $f_1(x)=x$ and $f_2(x)=x^2$ by functions in $C_c^\infty(\rr)$ coinciding with $f_1$, $f_2$ on $[-\widetilde{A},\widetilde{A}]$ for increasing values of $\widetilde{A}\in\nn$ and applying Proposition 4.6 in chapter 5 of \cite{ks} we conclude that there exist probability spaces on which processes of the same law as $Z_1$ and $Z_2$ (which we will also denote by $Z_1$ and $Z_2$) are defined such that
\begin{eqnarray}
dZ_i(t)=\mu(F_{\nu_i(t)}(Z_i(t)))\;dt+\sigma(F_{\nu_i(t)}(Z_i(t)))\;dW_i(t),\quad t\in[0,T]
\end{eqnarray}
holds for $i\in\{1,2\}$ and appropriate standard Brownian motions $W_1$, $W_2$. Next, we fix arbitrary numbers $x_1<x_2$ in $\rr$ and $t_1<t_2$ in $[0,T]$ and introduce the function $f(t,x)=\int_x^\infty \psi(t,y)\;dy$ on $[t_1,t_2]\times\rr$ where $\psi$ is an arbitrary continuous function on $[t_1,t_2]\times\rr$ which is continuously differentiable in both variables with $\psi(.,x)=\frac{\partial\psi}{\partial x}(.,x)=0$ whenever $x\notin(x_1,x_2)$. Applying Ito's formula to $f(t,Z_i(t))$ and taking the expectation we obtain
\begin{eqnarray*}
&&(\nu_i(t_2),f(t_2,.))-(\nu_i(t_1),f(t_1,.))\\
&&=\int_{t_1}^{t_2} \Big(\nu_i(t),\frac{\partial f}{\partial t}(t,.)+\frac{\partial f}{\partial x}(t,.)\mu(F_{\mathcal \nu_i(t)}(.))
+\frac{1}{2}\frac{\partial^2 f}{\partial x^2}(t,.)\sigma(F_{\mathcal \nu_i(t)}(.))^2\Big)\;dt
\end{eqnarray*}
for $i\in\{1,2\}$. Recalling that $\nu_1$ and $\nu_2$ are solutions of the Fokker-Planck equations in step 1 and following the proof of Lemma A.2 in \cite{ga} we conclude that the finite measures corresponding to the functionals $h\mapsto\int_0^T\int_\rr h(t,x)\;\nu_1(t)(dx)\;dt$, $h\mapsto\int_0^T\int_\rr h(t,x)\;\nu_2(t)(dx)\;dt$ acting on continuous bounded functions on $[0,T]\times\rr$ are absolutely continuous with respect to the Lebesgue measure on $[0,T]\times\rr$. Moreover, the proof of Lemma A.2 in \cite{ga} shows that the corresponding density functions $k_1$, $k_2$ on $[0,T]\times\rr$ are locally square integrable. Setting $w_i(t,x)=F_{\nu_i(t)}(x)$ for $t\in[0,T]$, $x\in\rr$, $i\in\{1,2\}$ and applying integration by parts with respect to the spatial variable in the last equation (recalling from step 1 that the measures $\nu_i(t)$, $t\in[0,T]$, $i\in\{1,2\}$ are absolutely continuous with respect to the Lebesgue measure on $\rr$) we see that for $i\in\{1,2\}$ it holds 
\begin{eqnarray*}
&&\int_\rr \psi(t_2,x) w_i(t_2,x)\;dx-\int_\rr \psi(t_1,x) w_i(t_1,x)\;dx\\
&&=\int_{t_1}^{t_2} \int_\rr \Big(\frac{\partial\psi}{\partial t}(t,x)w_i(t,x)+\frac{\partial\psi}{\partial x}(t,x)\Theta(w_i(t,x))\Big)\;dx\;dt\\
&&-\int_{t_1}^{t_2} \int_\rr \frac{\partial\psi}{\partial x}(t,x)\Sigma'(w_i(t,x))k_i(t,x)\;dx\;dt
\end{eqnarray*}
which we will call equation (*). Hereby, $\Theta$ and $\Sigma$ are the antiderivatives of $\mu$ and $\frac{1}{2}\sigma^2$, respectively, for which $\Theta(0)=\Sigma(0)=0$. Next, we note that for any $n\in\nn$ and any continuous function $\widetilde{\psi}:[t_1,t_2]\times\rr\rightarrow\rr$ which is supported on $D=[t_1,t_2]\times[x_1,x_2]$ and is continuously differentiable in $t$ and twice continuously differentiable in $x$ on $D$ we can find a function $\psi_n$ of the same type as the function $\psi$ in equation (*) such that $\psi_n\rightarrow_{n\rightarrow\infty}\psi$ uniformly on $D$ and 
\begin{eqnarray*}
\Big\|\psi_n-\widetilde{\psi}\Big\|_{L^2(D)}
+\Big\|\frac{\partial\psi_n}{\partial t}-\frac{\partial\widetilde{\psi}}{\partial t}\Big\|_{L^2(D)}
+\Big\|\frac{\partial\psi_n}{\partial x}-\frac{\partial\widetilde{\psi}}{\partial x}\Big\|_{L^2(D)}<\frac{1}{n}.
\end{eqnarray*}
This can be achieved by modifying $\widetilde{\psi}$ on small neighborhoods of 
\begin{eqnarray*}
([t_1,t_2]\times\{x_1\})\cup([t_1,t_2]\times\{x_2\}) 
\end{eqnarray*}
in $D$. Thus, an approximation argument together with the Cauchy-Schwarz inequality shows that equation (*) holds for all functions $\widetilde{\psi}$ of the described type and any $i\in\{1,2\}$. Choosing $k_i$ such that for each $t\in[0,T]$ the function $k_i(t,.)$ is a density function of the probability measure $\nu_i(t)$ for $i\in\{1,2\}$ and applying integration by parts to the last term in equation (*) we see that $w_i$, $i\in\{1,2\}$ are generalized solutions of the Cauchy problem for the generalized porous medium equation with non-linearity $\Sigma$ and convection $-\Theta$ in the sense of Definition 4 in \cite{gi}. From Theorem 4 in \cite{gi} we deduce $w_1=w_2$ by noting that $\mu$ and $\sigma$ can be easily extended to the whole of $\rr_+$ without violating Hypothesis 1 of \cite{gi}. Hence, $\nu_1(t)=\nu_2(t)$ for all $t\in[0,T]$ as desired. \ep

\section*{Acknowledgements}

The author thanks Amir Dembo for his invaluable comments throughout the preparation of this work. He is also grateful to Lenya Ryzhik and Andras Vasy for their suggestions.

\end{document}